\documentclass[10pt]{article}
\usepackage{amsmath,amssymb,amsfonts}
\newtheorem{theorem}{Theorem}
\newtheorem{lemma}{Lemma}
\newtheorem{proposition}{Proposition}
\renewcommand{\theequation}{\arabic{equation}}

\date{}

\title{Asymptotic Analysis of Expectations of Plane Partition Statistics}
\bigskip

\author{{\bf Ljuben Mutafchiev}\\
American University in Bulgaria, 2700 Blagoevgrad, Bulgaria \\ and
Institute of Mathematics and Informatics of the \\ Bulgarian
Academy of Sciences
\\ \tt {ljuben@aubg.bg}}

\begin{document}
\maketitle

\begin{abstract}
 Assuming that a plane partition of the positive integer $n$ is
chosen uniformly at random from the set of all such partitions, we
propose a general asymptotic scheme for the computation of
expectations of various plane partition statistics as $n$ becomes
large. The generating functions that arise in this study are of
the form $Q(x)F(x)$, where $Q(x)=\prod_{j=1}^\infty (1-x^j)^{-j}$
is the generating function for the number of plane partitions. We
show how asymptotics of such expectations can be obtained directly
from the asymptotic expansion of the function $F(x)$ around $x=1$.
The representation of a plane partition as a solid diagram of
volume $n$ allows interpretations of these statistics in terms of
its dimensions and shape. As an application of our main result, we
obtain the asymptotic behavior of the expected values of the
largest part, the number of columns, the number of rows (that is,
the three dimensions of the solid diagram) and the trace (the
number of cubes in the wall on the main diagonal of the solid
diagram). Our results are similar to those of Grabner et al.
\cite{GKW14} related to linear integer partition statistics. We
base our study on the Hayman's method for admissible power series.
\end{abstract}

\vspace{.5 cm}



 \vspace{.2cm}

\section{Introduction}

A plane partition $\omega$ of the positive integer $n$ is an array
of non-negative integers
\begin{equation}\label{planepart}
\begin{array}{cccccccccccc}
\omega_{1,1} & \omega_{1,2} & \omega_{1,3} \quad \cdots \\
\omega_{2,1} & \omega_{2,2} & \omega_{2,3} \quad \cdots \\
\cdots & \cdots & \cdots \\
\end{array}
\end{equation}
that satisfy $\sum_{h,j\ge 1}\omega_{h,j}=n$,
and the rows and columns in (\ref{planepart}) are arranged in decreasing order:
$\omega_{h,j}\ge\omega_{h+1,j}$ and $\omega_{h,j}\ge\omega_{h,j+1}$ for all $h,j\ge 1$.
The non-zero entries $\omega_{h,j}>0$ are called parts of $\omega$. If there are
$\lambda_h$ parts in the $h$th row, so that for some $l$,
$\lambda_1\ge\lambda_2\ge...\ge\lambda_l >\lambda_{l+1}=0$, then the (linear)
partition $\lambda=(\lambda_1,\lambda_2,...,\lambda_l)$ of the integer
$m=\lambda_1+\lambda_2+...+\lambda_l$ is called the shape of $\omega$, denoted by
$\lambda$. We also say that $\omega$ has $l$ rows and $m$ parts. Sometimes, for the sake of brevity,
the zeroes in array (\ref{planepart}) are deleted. For instance, the abbreviation
$$
\begin{array}{cccc}
5 & 4 & 1 & 1 \\
3 & 2 & 1 \\
2 & 1 \\
\end{array}
$$
represents a plane partition of $n=20$ with $l=3$ rows and $m=9$ parts.

Any plane partition $\omega$ has an associated solid diagram
$\Delta=\Delta(\omega)$ of volume $n$. It is defined as a set of
$n$ integer lattice points
$\mathbf{x}=(x_1,x_2,x_3)\in\mathbb{N}^3$, such that if
$\mathbf{x}\in\Delta$ and $x_j^\prime\le x_j, j=1,2,3$, then
$\mathbf{x}^\prime=(x_1^\prime,x_2^\prime,x_3^\prime)\in\Delta$ too.
(Here $\mathbb{N}$ denotes the set of all positive integers.)
Indeed the entry $\omega_{h,j}$ can be interpreted as the height of the column of unit cubes
stacked along the vertical line $x_1=h, x_2=j$, and the solid diagram is the union of
all such columns.

Plane partitions are originally introduced by Young \cite{Y01} as
a natural generalization of integer partitions in the plane. The
problem of enumerating plane partitions was studied first by
MacMahon \cite{M12} (see also \cite{M16}), who showed that, for
any parallelepiped $B(l,s,t)=\{(h,j,k)\in\mathbb{N}^3: h\le l,
j\le s, k\le t\}$ and any $\mid x\mid<1$,
\begin{equation}\label{mcmahon}
\sum_{\Delta\subset B(l,s,t)} x^{\mid\Delta\mid} =\prod_{(h,j,k)\in B(l,s,t)}
\frac{1-x^{h+j+k-1}}{1-x^{h+j+k-2}},
\end{equation}
where $\Delta$ is the solid diagram defined above and
$\mid\Delta\mid$ denotes its volume (a more explicit version of
(\ref{mcmahon}) may be found in \cite{P05}).

Let $q(n)$ denote the total number of plane partitions of the positive integer $n$
(or, the total number of solid diagrams of volume $n$). It turns out that (\ref{mcmahon})
implies the following generating function identity:
\begin{equation}\label{gfq}
Q(x)=1+\sum_{n=1}^\infty q(n)x^n =\prod_{j=1}^\infty(1-x^j)^{-j}
\end{equation}
(more details may be also found, e.g., in
\cite[Corollary~18.2]{S71} and \cite[Corollary~11.3]{A76}).
Subsequent research on the enumeration of plane partitions was
focused on bijective interpretations and proofs of MacMahon's
formula (\ref{gfq}) (see, e.g., \cite{K99,P01}).

The asymptotic form of the numbers $q(n)$, as $n\to\infty$, has
been obtained by Wright \cite{W31} (see also \cite{MK06} for a
little correction). It is given by the following formula:
\begin{equation}\label{wright}
q(n)\sim \frac{(\zeta(3))^{7/36}}{2^{11/36}(3\pi)^{1/2}} n^{-25/36}
\exp{(3(\zeta(3))^{1/3}(n/2)^{2/3}+2\gamma)},
\end{equation}
where
$$
\zeta(z)=\sum_{j=1}^\infty j^{-z}
$$
is the Riemann zeta function and
\begin{equation}\label{glakin}
\gamma=\int_0^\infty\frac{u\log{u}}{e^{2\pi u}-1}du
=\frac{1}{2}\zeta^\prime(-1).
\end{equation}
(The constant $\zeta^\prime(-1)=-0.1654...$ is closely related to
Glaisher-Kinkelin constant; see \cite{F03}).

{\it Remark 1.} In fact, Wright \cite{W31} has obtained an
asymptotic expansion for $q(n)$ using the circle method.

Next, we introduce the uniform probability measure $\mathbb{P}$ on
the set of plane partitions of $n$, assuming that the probability
$1/q(n)$ is assigned to each plane partition. In this way, each
numerical characteristic of a plane partition of $n$ becomes a
random variable (a statistic in the sense of the random generation
of plane partition of $n$). In the following, we will discuss
several different instances of plane partition statistics. Our
goal is to develop a general asymptotic scheme that allows us to
derive an asymptotic formula for the $n$th coefficient
$[x^n]Q(x)F(x)$ of the product $Q(x)F(x)$, where $Q(x)$ is defined
by (\ref{gfq}) and the power series $F(x)$ is suitably restricted
on its behavior in a neighborhood of $x=1$. We will show further
that expectations of plane partitions statistics we will consider
lead to generating functions of this form. From one side, our
study is motivated by the asymptotic results of Grabner et al.
\cite{GKW14} on linear partitions statistics. Their study is based
on general asymptotic formulae for the $n$th coefficient of a
similar product of generating functions with $Q(x)$ replaced by
the Euler partition generating function $P(x)=\prod_{j=1}^\infty
(1-x^j)^{-1}.$ The second factor $F(x)$ satisfies similar general
analytic conditions around $x=1$. In addition, our interest to
study plane partitions statistics was also attracted by several
investigations during the last two decades on shape parameters of
random solid diagrams of volume $n$ as $n\to\infty$. Below we
present a brief account on this subject.

Cerf and Kenyon \cite{CK01} have determined the asymptotic shape
of the random solid diagram, while Cohn et al. \cite{CLP98} have
studied a similar problem whenever a solid diagram is chosen
uniformly at random among all diagrams boxed in $B(l,s,t)$, for
large $l, s$ and $t$, all of the same order of magnitude. Okounkov
and Reshetikhin \cite{OR03} rediscovered Cerf and Kenyon's
limiting shape result and studied asymptotic correlations in the
bulk of the random solid diagram. Their analysis is based on a
deterministic formula for the correlation functions of the Schur
process. The joint limiting distribution of the height (the
largest part in (\ref{planepart})), depth (the number of columns
in (\ref{planepart})) and width (the number of rows in
(\ref{planepart})) in a random solid diagram was obtained by
Pittel \cite{P05}. The one-dimensional marginal limiting
distributions of this random vector were found in \cite{M06}. The
trace of a plane partition is defined as the sum of the diagonal
parts in (\ref{planepart}). Its limiting distribution is
determined in \cite{KM07}. Bodini et al. \cite{BFP10} studied
random generators of plane partitions according to the size of
their solid diagrams. They obtained random samplers that are of
complexity $O(n\log^3{n})$ in an approximate-size sampling and of
complexity $O(n^{4/3})$ in exact-size sampling. These random
samplers  allow to perform simulations in order to confirm the
known results about the limiting shape of the plane partitions.

In the proof of our main asymptotic result for the coefficient
$[x^n]Q(x)F(x)$ we use the saddle point method. In contrast to
\cite{GKW14}, we base our study on a theorem due to Hayman
\cite{H56} for estimating coefficients of admissible power series
(see also \cite[Section~VIII.5]{FS09}). We show that $Q(x)$ (see
(\ref{gfq})) is Hayman admissible function and impose conditions
on $F(x)$ which are given in terms of the Hayman's theorem. In the
examples we present, we demonstrate two different but classical
approaches for estimating power series around their main
singularity.

Our paper is organized as follows. In Section 2 we include the
Hayman's admissibility conditions and his main asymptotic result.
We briefly describe a relationship between Hayman's theorem and
Meinardus approach \cite{M54} for obtaining the asymptotic
behavior of the Taylor coefficients of infinite products of the
form
\begin{equation}\label{prodbj}
f_b(x)=\prod_{j=1}^\infty (1-x^j)^{-b_j}
\end{equation}
under certain general conditions on the sequence $\{b_j\}_{j\ge
1}$ of non-negative numbers. In Section 2, we also state our main
result (Theorem 1) for the asymptotic of the $n$th coefficient
$[x^n]Q(x)F(x)$ under certain relatively mild conditions on $F(x)$. The
proof of Theorem 1 is given in Section 3. Section 4 contains some
examples of plane partition statistics that lead to generating
functions of the form $Q(x)F(x)$. We apply Theorem 1 to obtain the
asymptotic behavior of expectations of the underlying statistics.
For the sake of completeness, in the Appendix we show how Wright's
formula (\ref{wright}) follows from Hayman's theorem.

\section{Some Remarks on Hayman Admissible Functions and Meinardus Theorem on Weighted Partitions. Statement of the Main Tesult}

Our first goal in this section is to give a brief introduction to
the analytic combinatorics background that we will use in the
proof of our main result. Clearly, $x^n[Q(x)F(x)]$ with $Q(x)$
given by (\ref{gfq}) can be represented by a Cauchy integral whose
integrand includes the product $Q(x)F(x)$ (the conditions that
$F(x)$ should satisfy will be specified later). Its asymptotic
behavior heavily depends on the analytic properties of $Q(x)$
whose infinite product representation (\ref{gfq}) shows that the
unit circle is a natural boundary and its main singularity is at
$x=1$. The main tools for the asymptotic analysis of
$q(n)=x^n[Q(x)]$ are either based on the circle method (see
\cite{W31}) or on the saddle-point method (see \cite{M54} and
\cite{BFP10}). Both yield Wright's asymptotic formula
(\ref{wright}). An asymptotic formula in a more general framework
(see (\ref{prodbj})) was obtained by Meinardus \cite{M54} (see
also \cite{GSE08} for some extensions). A proof of formula
(\ref{wright}) that combines Meinardus approach with Hayman's
theorem for admissible power series is started in Section 3 and
completed in the Appendix.

Here we briefly describe Menardus' approach, which is essentially
based on analytic properties of the Dirichlet generating series
$$
D(z)=\sum_{j=1}^\infty b_j j^{-z}, \quad z=u+iv,
$$
where the sequence of non-negative numbers $\{b_j\}_{j\ge 1}$ is
the same as in the infinite product (\ref{prodbj}). We will avoid
the precise statement of Meinardus' assumptions on $\{b_j\}_{j\ge
1}$ as well as some extra notations and concepts. The first
assumption $(M_1)$ specifies the domain $\mathcal{H}=\{z=u+iv:
u\ge -C_0, 0<C_0<1\}$ in the complex plane, in which $D(z)$ has an
analytic continuation. The second one $(M_2)$ is related to the
asymptotic behavior of $D(z)$ whenever $\mid v\mid\to\infty$. A
function of the complex variable $z$ which is bounded by
$O(\mid\Im{(z)}\mid^{C_1}), 0<C_1<\infty$, in certain domain in
the complex plane is called function of finite order. Meinardus'
second condition $(M_2)$ requires that $D(z)$ is of finite order
in the whole domain $\mathcal{H}$. Finally, the Meinardus' third
condition $(M_3)$ implies a bound on the ordinary generating
function of the sequence $\{b_j\}_{j\ge 1}$. It can be stated in a
way, simpler than the Meinardus' original expression, by the
inequality:
$$
\sum_{j=1}^\infty b_j e^{-j\alpha}\sin^2{(\pi ju)}\ge C_2
\alpha^{-\epsilon_0}, \quad 0<\frac{\alpha}{2\pi}<\mid
u\mid<\frac{1}{2},
$$
for sufficiently small $\alpha$ and some constants $C_2,
\epsilon_0>0$ $(C_2=C_2(\epsilon_0))$ (see \cite[p.~310]{GSE08}).

The infinite product representation (\ref{gfq}) for $Q(x)$ implies
that $b_j=j, j\ge 1$, and therefore, $D(z)=\zeta(z-1)$. It is
known that this sequence satisfies the Meinardus scheme of
conditions (see, e.g., \cite{MK06} and \cite[p.~312]{GSE08}).

Now we proceed to Hayman admissibility method \cite{H56},
\cite[Section~VIII.5]{FS09}. To present the idea and show how it
can be applied to the proof of our main result, we need to
introduce some auxiliary notations.

Consider a function $G(x)=\sum_{n=0}^\infty g_n x^n$ that is
analytic for $\mid x\mid<\rho, 0<\rho\le\infty$. For $0<r<\rho$,
we set
\begin{equation}\label{a}
a(r)=r\frac{G^\prime(r)}{G(r)},
\end{equation}
\begin{equation}\label{b}
b(r)=\frac{rG^\prime(r)}{G(r)}+r^2\frac{G^{\prime\prime}(r)}{G(r)}
-r^2\left(\frac{G^\prime(r)}{G(r)}\right)^2.
\end{equation}
In the statement of Hayman's result we use the terminology given
in \cite[Section~VIII.5]{FS09}. We assume that $G(x)>0$ for $x\in
(R_0,\rho)\subset (0,\rho)$ and satisfies the following three
conditions:

{\it Capture condition.} $lim_{r\to\rho} a(r)=\infty$ and
$\lim_{r\to\rho} b(r)=\infty$.

{\it Locality condition.} For some function $\delta=\delta(r)$
defined over $(R_0,\rho)$ and satisfying $0<\delta<\pi$, one has
$$
G(r e^{i\theta})\sim G(r) e^{i\theta a(r)-\theta^2 b(r)/2}
$$
as $r\to\rho$, uniformly for $\mid\theta\mid\le\delta(r)$.

{\it Decay condition.}
$$
G(r e^{i\theta}) =o\left(\frac{G(r)}{\sqrt{b(r)}}\right)
$$
as $r\to\rho$ uniformly for $\delta(r)<\mid\theta\mid\le\pi$.

{\bf Hayman Theorem.} Let $G(x)$ be Hayman admissible function and
$r=r_n$ be the unique solution of the equation
\begin{equation}\label{arn}
a(r)=n.
\end{equation}
Then the Taylor coefficients $g_n$ of $G(x)$ satisfy, as
$n\to\infty$,
\begin{equation}\label{gn}
g_n\sim\frac{G(r_n)}{r_n^n\sqrt{2\pi b(r_n)}}
\end{equation}
with $b(r_n)$ given by (\ref{b}).

In the next section we will show that $Q(x)$, given by
(\ref{gfq}), is admissible in the sense of Hayman and apply Hayman
Theorem setting $G(x):=Q(x)$.

Now, we proceed to the statement of our main result. As in
\cite{GKW14}, we assume that $F(x)$ satisfies two rather mild
conditions. Since we will employ Hayman's method, the first one is
given in terms of eq. (\ref{arn}). The second one requires, as in
\cite{GKW14}, that $F(x)$ does not grow too fast as $\mid x\mid\to
1$.

{\it Condition A.} Let $r=r_n$ be the solution of (\ref{arn}). We
assume that
$$
\lim_{n\to\infty}\frac{F(r_n e^{i\theta})}{F(r_n)}=1
$$
uniformly for $\mid\theta\mid\le\delta(r_n)$, where $\delta(r)$ is
the function defined by Hayman's locality condition.

{\it Condition B.} There exist two constants $C>0$ and $\eta\in
(0,2/3)$, such that, as $\mid x\mid\to 1$,
$$
F(x)=O(e^{C/(1-\mid x\mid)^\eta}).
$$

Our main result is as follows.

\begin{theorem} Let $\{d_n\}_{n\ge 1}$ be a sequence with the following expansion:
\begin{equation}\label{dn}
d_n=\left(\frac{2\zeta(3)}{n}\right)^{1/3} -\frac{1}{36n}
+O(n^{-1-\beta}), \quad n\to\infty,
\end{equation}
where $\beta>0$ is certain fixed constant. Furthermore, suppose
that the function $F(x)$ satisfies conditions (A) and (B) and
$Q(x)$ is the infinite product  given by (\ref{gfq}). Then, there
is a constant $c>0$ such that, as $n\to\infty$,
$$
\frac{1}{q(n)}[x^n]Q(x)F(x)=F(e^{-d_n})(1+o(1))
+O(e^{-cn^{2/9}/\log^2{n}}),
$$
where $q(n)$ is the $n$th coefficient  in the Taylor expansion
of $Q(x)$.
\end{theorem}

The proof that we will present in the next section stems from the
Cauchy coefficient formula in which the contour of integration is
the circle $x=e^{-d_n+i\theta}, -\pi<\theta\le\pi$ and $d_n$ is
defined by (\ref{dn}). We have
\begin{equation}\label{cauchy}
[x^n]Q(x)F(x) =\frac{e^{nd_n}}{2\pi}\int_{-\pi}^\pi
Q(e^{-d_n+i\theta})F(e^{-d_n+i\theta}) e^{-i\theta n} d\theta.
\end{equation}
The proof of Theorem 1 is divided into two parts:

(i) Proof of Hayman admissibility for $Q(x)$.

(ii) Obtaining an asymptotic estimate for the Cauchy integral
(\ref{cauchy}).

\section{Proof of Theorem 1}

{\it Part (i).}

We will essentially use some more general observations established
in \cite{GSE08, M13}.

First, we set in (\ref{a}) and (\ref{b}) $G(x):=Q(x)$ and
$r=r_n:=e^{-d_n}$. The next lemma is a particular case of a more
general result due to Granovsky et al. \cite[Lemma~2]{GSE08}.

\begin{lemma} For enough large $n$, the unique solution of the equation
$$
a(e^{-d_n})=n
$$
is given by (\ref{dn}). Moreover, as $n\to\infty$,
\begin{equation}\label{bdn}
b(e^{-d_n})\sim\frac{3n^{4/3}}{(2\zeta(3))^{1/3}}.
\end{equation}
\end{lemma}

{\it Proof.} The Dirichlet generating series of the sequence
$\{j\}_{j\ge 1}$ is $\zeta(z-1)$, which is a meromorphic function
in the complex plane with a single pole at $z=2$ with residue $1$
(see, e.g., \cite[Section~13.13]{WW27}). Granovsky et al.
\cite{GSE08} showed thay it satisfies Meinardus' conditions
$(M_1)$ and $(M_2)$. Hence, by formula (43) in \cite{GSE08}, we
have
\begin{eqnarray}
& & d_n =(\Gamma(2)\zeta(3))^{1/3} n^{-1/3} +\frac{\zeta(-1)}{3n}
+O(n^{-1-\beta}) \nonumber \\
& & =\left(\frac{2\zeta(3)}{n}\right)^{1/3} -\frac{1}{36n}
+O(n^{-1-\beta}), \quad \beta>0, \nonumber
\end{eqnarray}
which is just (\ref{dn}) (here we have also used that
$\zeta(-1)=-1/12$ \cite[Section~13.14]{WW27}). Formula (\ref{bdn})
follows from (\ref{dn}) and a more general result from
\cite[formula~(2.6)]{M13}. $\Box$

Lemma 1 shows that $a(e^{-d_n})\to\infty$ and $b(e^{-d_n})\to\infty$
as $n\to\infty$, that is Hayman's ''capture'' condition is satisfied
with $r=r_n=e^{-d_n}$. To show next that Hayman's ''decay'' condition
is satisfied by $Q(x)$, we set
\begin{equation}\label{delta}
\delta_n=\frac{d_n^{5/3}}{\log{n}} =\frac{1}{\log{n}}\left(\frac{2\zeta(3)}{n}\right)^{5/9}
(1+O(n^{-2/3}))
\end{equation}
with $d_n$ given by (\ref{dn}). The next lemma is a particular
case of Lemma 2.4 in \cite{M13}.

\begin{lemma}
For sufficiently large $n$, we have
$$
\mid Q(e^{-d_n+i\theta}\mid\le Q(e^{-d_n})e^{-C_3 d_n^{-2/3}}
$$
uniformly for $\delta_n\le\mid\theta\mid<\pi$, where $C_3>0$ is
an absolute constant and
$d_n$ is defined by $(\ref{dn})$.
\end{lemma}

{\it Proof.} We will apply Lemma 2.4 from \cite{M13} with
$f_b(x):=Q(x)$ (see (\ref{prodbj})) and and $\alpha_n:=d_n$. It
shows that the ratio $\mid f_b(e^{-d_n+i\theta})\mid/f(e^{-d_n})$
is bounded by $e^{-C_3 d_n^{-\epsilon_1}}$ uniformly for
$\delta_n\le\mid\theta\mid<\pi$, where $C_3, \epsilon_1>0$. Our
goal is to show that the plane partition generating function
$Q(x)$ is bounded in the same way with $\epsilon_1=2/3$. To show
this, we notice first that the chain of inequalities (2.18) in
\cite{M13} with $b_j=j, j\ge 1$ implies that
\begin{equation}\label{relog}
\Re{(\log{Q(e^{-d_n+i\theta})})} -\log{Q(e^{-d_n})}\le\frac{\log{5}}{2} S_n
\end{equation}
uniformly for $d_n\le\mid\theta\mid<\pi$, where
$$
S_n=\sum_{j=1}^\infty je^{-j d_n}\sin^2{(\pi j\theta/2)}.
$$
Recall that we deal with the Dirichlet generating function
$\zeta(z-1)$ that satisfies the conditions of Lemma 1 in
\cite{GSE08} (namely, it converges in the half-plane
$\Re{(z)}>2>1$). So, we can apply its last part in combination
with (\ref{relog}) to conclude that there exists a constant
$C_3^\prime>0$ such that
\begin{equation}\label{gran}
\frac{\mid Q(e^{-d_n+i\theta})\mid}{Q(e^{-d_n})}
=\exp{(\Re{(\log{Q(e^{-d_n+i\theta})}))} -\log{Q(e^{-d_n})}} \le
e^{-C_3^\prime d_n^{-1}}
\end{equation}
uniformly for $d_n\le\mid\theta\mid<\pi$. Furthermore, since the
abscissa of convergence of $\zeta(z-1)$ is $2$ (and $2+1<\pi$),
for $\delta_n\le\mid\theta\mid<d_n$, we can apply the estimate
(2.21) from \cite[p.~440]{M13} with $\omega(n)=\log{n}$.
Therefore, for enough large $n$, we have
$$
S_n\ge C_3^{\prime\prime}d_n^{-2/3}/\log^2{n},
$$
where $C_3^{\prime\prime}>0$. Hence by (\ref{relog}), uniformly
for $\delta_n\le\mid\theta\mid<d_n$,
\begin{equation}\label{grann}
\frac{\mid Q(e^{-d_n+i\theta})\mid}{Q(e^{-d_n})} \le
e^{-C_3^{\prime\prime}d_n^{-2/3}/\log^2{n}}.
\end{equation}
Combining (\ref{gran}) and (\ref{grann}) and setting
$C_3=\min{(C_3^\prime,C_3^{\prime\prime})}$, we obtain the
required uniform estimate for $\delta_n\le\mid\theta\mid<\pi$.
$\Box$

This lemma, in combination with (\ref{bdn}) and (\ref{dn}),
implies that $\mid
Q(e^{-d_n+i\theta})\mid=o(Q(e^{-d_n})/\sqrt{b(e^{-d_n})})$
uniformly for $\delta_n\le\mid\theta\mid\le\pi$, which is just
Hayman's ''decay'' condition.

Finally, since $\zeta(z-1)$ satisfies Meinardus' conditions
$(M_1)$ and $(M_2)$, Lemma 2.3 from \cite{M13}, implies Haiman's
''locality'' condition for $Q(x)$.

\begin{lemma} We have,
$$
e^{-i\theta n}\frac{Q(e^{-d_n+i\theta})}{Q(e^{-d_n})}
=e^{-\theta^2 b(e^{-d_n})/2}(1+O(1/\log^3{n}))
$$
as $n\to\infty$ uniformly for $\mid\theta\mid\le\delta_n$, where
$\delta_n$, $d_n$ and $b(e^{-d_n})$ are determined by
(\ref{delta}), (\ref{dn}) and (\ref{b}), respectively.
\end{lemma}

So, all conditions of Hayman's theorem hold, and we can apply it
with $g_n:=q(n), G(x):=Q(x), r_n:=e^{-d_n}$ and $\rho=2$ to find
that
\begin{equation}\label{qhay}
q(n)\sim\frac{e^{n d_n}Q(e^{-d_n})}{\sqrt{2\pi b(e^{-d_n})}},
\quad n\to\infty.
\end{equation}
In the Appendix we will show that (\ref{qhay}) implies the
corrected version of Wright's formula (\ref{wright}).

{\it Remark 2.} The fact that the MacMahon's generating function
$Q(x)$ given by (\ref{gfq}) is admissible in the sense of Hayman
is a particular case of a more general result established in
\cite{M13}  and related to the infinite products $f_b(x)$ of the
form (\ref{prodbj}). It turns out that the Meinardus'scheme of
assumptions on $\{b_j\}_{j\ge 1}$ implies that $f_b(x)$ is
admissible in the sense of Hayman.

{\it Part (ii)}

We break up the range of integration in (\ref{cauchy}) as follows:
\begin{equation}\label{twiint}
[x^n]Q(x)F(x) =J_{1,n} +J_{2,n},
\end{equation}
where
\begin{equation}\label{jone}
J_{1,n} =\frac{e^{nd_n}}{2\pi} \int_{-\delta_n}^{\delta_n}
Q(e^{-d_n+i\theta})F(e^{-d_n+i\theta}) e^{-i\theta n} d\theta,
\end{equation}
\begin{equation}\label{jtwo}
J_{2,n} =\int_{\delta_n\le\mid\theta\mid<\pi}
Q(e^{-d_n+i\theta})F(e^{-d_n+i\theta}) e^{-i\theta n} d\theta
\end{equation}
and $\delta_n$ is defined by (\ref{delta}).

The estimate for $J_{1,n}$ follows from Hayman's ''locality''
condition and condition (A) for $F(x)$. We have
\begin{eqnarray}\label{joneasymp}
& & J_{1,n} =\frac{e^{n d_n}Q(e^{-d_n})F(e^{-d_n})}{2\pi}
\int_{-\delta_n}^{\delta_n}
\left(\frac{Q(e^{-d_n+i\theta})}{Q(e^{-d_n})}\right) e^{-i\theta
n} \left(\frac{F(e^{-d_n+i\theta})}{F(e^{-d_n})}\right) d\theta
\nonumber \\
& & =\frac{e^{n d_n}Q(e^{-d_n}) F(e^{-d_n})}{2\pi}
\int_{-\delta_n}^{\delta_n} e^{-\theta^2 b(e^{-d_n})/2}
\left(1+O\left(\frac{1}{\log^3{n}}\right)\right)(1+o(1)) d\theta
\nonumber \\
& & \sim \frac{e^{n d_n}Q(e^{-d_n}) F(e^{-d_n})}{2\pi}
\int_{-\delta_n}^{\delta_n} e^{-\theta^2 b(e^{-d_n})/2} d\theta.
\end{eqnarray}
Substituting $\theta=y/\sqrt{b(e^{-d_n})}$, we observe that
\begin{eqnarray}
& & \int_{-\delta_n}^{\delta_n} e^{-\theta^2 b(e^{-d_n})/2}
d\theta \sim\frac{1}{\sqrt{b(e^{-d_n})}} \int_{-\delta_n
\sqrt{b(e^{-d_n})}}^{\delta_n\sqrt{b(e^{-d_n})}} e^{-y^2/2} dy
\nonumber \\
 & & \sim\int_{-\infty}^\infty
e^{-y^2/2}dy=\sqrt{\frac{2\pi}{b(e^{-d_n})}}, \quad n\to\infty,
\nonumber
\end{eqnarray}
since by (\ref{bdn}) and (\ref{delta})
$$
\delta_n\sqrt{b(e^{-d_n})}\sim\sqrt{3}(2\zeta(3))^{7/18}
\frac{n^{1/9}}{\log{n}}, \quad n\to\infty.
$$
Hence, by (\ref{qhay}), the asymptotic estimate (\ref{joneasymp})
is simplified to
\begin{equation}\label{jonefin}
J_{1,n} =q(n)F(e^{-d_n}) +o(q(n)F(e^{-d_n})).
\end{equation}

To estimate $J_{2,n}$, we first apply Lemma 2 with $d_n$ replaced
by its expression (\ref{dn}). Thus we observe that there is a
constant $C_4>0$ such that the inequality
\begin{equation}\label{ineqq}
\mid Q(e^{-d_n+i\theta})\mid\le Q(e^{-d_n}) e^{-C_4
n^{2/9}/\log^2{n}}
\end{equation}
holds uniformly for $\delta_n\le\mid\theta\mid<\pi$. We will
combine this estimate with condition (B) on the function $F(x)$.
It implies that, for certain constants $c_0, c_1>0$, we have
\begin{equation}\label{efest}
F(e^{-d_n})=O(e^{c_0 d_n^{-\eta}}) =O(e^{c_1 n^{\eta/3}}).
\end{equation}
Now, combining (\ref{jtwo}), (\ref{ineqq}), (\ref{bdn}),
(\ref{qhay}) and (\ref{efest}), we obtain
\begin{eqnarray}\label{jtwofin}
& & \mid J_{2,n}\mid\le \int_{\delta_n\le\mid\theta\mid<\pi} \mid
Q(e^{-d_n+i\theta})\mid \mid F(e^{-d_n+i\theta})\mid d\theta
\nonumber \\
& & \le\frac{e^{n d_n}}{\pi} Q(e^{-d_n})O(e^{c_1 n^{\eta/3}})
(\pi-\delta_n) e^{-C_4 n^{2/9}/\log{n}} \nonumber \\
& & =O\left(\frac{e^{n d_n} Q(e^{-d_n})}{\sqrt{2\pi b(e^{-d_n})}}
n^{2/3} \exp{\left(-\frac{C_4 n^{2/9}}{\log^2{n}} +c_1
n^{\eta/3}\right)}\right) \nonumber \\
& & =O(q(n) e^{-c n^{2/9}/\log^2{n}}),
\end{eqnarray}
for some $c>0$. Substituting the estimates obtained in
(\ref{jonefin}) and (\ref{jtwofin}) into (\ref{twiint}), we
complete the proof. $\Box$

\section{Examples}

{\bf The trace of a plane partition.} The trace $T_n$ of a plane
partition $\omega$, given by array (\ref{planepart}), is defined
as the sum of its diagonal parts:
$$
T_n =\sum_{j\ge 1} \omega_{j,j}.
$$
The asymptotic behavior of $T_n$, as $n\to\infty$, can be studied
using the following generating function identity established by
Stanley \cite{S73} (see also \cite[Chapter~11, Problem~5]{A76}):
\begin{equation}\label{bivgftr}
1+\sum_{n=1}^\infty q(n)x^n\sum_{m=1}^n \mathbb{P}(T_n=m)u^m
=1+\sum_{n=1}^\infty q(n)\varphi_n(u)x^n =\prod_{j=1}^\infty
(1-ux^j)^{-j},
\end{equation}
where $\varphi_n(u)$  denotes the probability generating function
of $T_n$: $\varphi_n(u) =\mathbb{E}(u^{T_n})$ ($\mid u\mid\le 1$).
Since $\varphi_n^\prime(1)=\mathbb{E}(T_n)$, a differentiation of
(\ref{bivgftr}) with respect to $u$ yields
\begin{equation}\label{expectfone}
\sum_{n=1}^\infty q(n)\mathbb{E}(T_n)x^n =Q(x)F_1(x),
\end{equation}
where
\begin{equation}\label{fone}
F_1(x) =\sum_{j=1}^\infty\frac{jx^j}{1-x^j}.
\end{equation}
For $\mid\theta\mid\le\delta(r_n)$, by Taylor formula we have
$$
F_1(r_n e^{i\theta}) =F_1(r_n) +O(\mid\theta\mid F_1^\prime(r_n))
=F_1(r_n)+O(\delta(r_n)F_1^\prime(r_n)),
$$
where $r_n$ is the solution of (\ref{arn}). So, the function
$F_1(x)$ satisfies condition (A) if
\begin{equation}\label{conda}
\frac{F_1^\prime(r_n)}{F_1(r_n)}\delta(r_n)\to 0, \quad
n\to\infty.
\end{equation}
Differentiating (\ref{expectfone}), we get
\begin{equation}\label{derivfone}
F_1^\prime(x)=\sum_{j=1}^\infty\frac{j^2 x^{j-1}}{(1-x^j)^2}.
\end{equation}
Setting in (\ref{fone}) and (\ref{expectfone}) $x=r_n$ and
interpreting the sums as Riemann sums with step size
$-\log{r_n}=-\log{(1-(1-r_n))}=1-r_n+O((1-r_n)^2)$, it is easy to
show that
$$
F_1^\prime(r_n)
=O\left((1-r_n)^{-3}\int_0^\infty\frac{u^2}{(e^u-1)^2}du \right)
=O((1-r_n)^{-3}).
$$
Hence, with $r_n=e^{-d_n}$, we have $1-r_n=d_n+O(d_n^2)$ and
\begin{equation}\label{derivdn}
F_1^\prime(x)\mid_{x=e^{-d_n}}=O(d_n^{-2}).
\end{equation}
For $F_1(e^{-d_n})$ we need a more precise estimate. In the same
way, using the Riemann sum approximation, we obtain
\begin{eqnarray}\label{fonedn}
& & F_1(e^{-d_n}) =\sum_{j=1}^\infty \frac{j e^{-j d_n}}{1-e^{-j
d_n}} =d_n^{-2}\sum_{j=1}^\infty
\frac{jd_n e^{-jd_n}}{1-e^{-jd_n}}d_n  \nonumber \\
& & \sim d_n^{-2}\int_0^\infty\frac{u}{e^u-1}du =d_n^{-2}\zeta(2),
\quad n\to\infty,
\end{eqnarray}
where in the last equality we have used formula 27.1.3 from
\cite{AS65}. Now, the convergence in (\ref{conda}) follows from
(\ref{derivdn}), (\ref{fonedn}) and (\ref{delta}) and thus
$F_1(x)$ satisfies condition (A). Condition (B) is also obviously
satisfied, since an argument similar to that in (\ref{fonedn})
implies that $F_1(x)=O((1-\mid x\mid)^{-2}) =O(e^{C/(1-\mid
x\mid)^\eta})$, as $\mid x\mid\to 1$, for any $C>0$ and $\eta\in
(0,2/3)$.

Combining (\ref{fonedn}) with (\ref{dn}) and applying the result
of Theorem 1 to (\ref{expectfone}) we obtain the following
asymptotic equivalence for $\mathbb{E}(T_n)$.

\begin{proposition} If $n\to\infty$, then
$$
\mathbb{E}(T_n)\sim\kappa_1 n^{2/3},
$$
where $\kappa_1=(2\zeta(3))^{-2/3}\pi^2/6=0.9166...$.
\end{proposition}

{\it Remark 3.} One can compare this asymptotic result with the
limit theorem for $T_n$ obtained in \cite{KM07}, where it is shown
that $T_n$, appropriately normalized, converges weakly to the
standard Gaussian distribution.

{\bf The largest part, the number of rows and the number of
columns of a plane partition.} Let $X_n$, $Y_n$ and $Z_n$ denote
the size of the largest part, the number of rows and number of
columns in a random plane partition of $n$, respectively. Using
the solid diagram interpretation $\Delta(\omega)$ of a plane
partition $\omega$, one can interpret $X_n$, $Y_n$ and $Z_n$ as
the height, width and depth of $\Delta(\omega)$. Any permutation
$\sigma$ of the coordinate axes in $\mathbb{N}^3$, different from
the identical one, transforms $\Delta(\omega)$ into a diagram that
uniquely determines another plane partition $\sigma\circ\omega$.
The permutation $\sigma$ also permutes the three statistics
$(X_n,Y_n,Z_n)$. So, if one of these statistics is restricted by
an inequality, the same restriction occurs on the statistic
permuted by $\sigma$. The one to one correspondence between
$\omega$ and $\sigma\circ\omega$ implies that $X_n$, $Y_n$ and
$Z_n$ are identically distributed for every fixed $n$ with respect
to the probability measure $\mathbb{P}$. (More details may be
found in \cite[p. 371]{S99}.) Hence, in the context of the
expected value $\mathbb{E}$ with respect to the probability
measure $\mathbb{P}$, we will use the common notation
$\mathbb{E}(W_n)$ for $W_n=X_n,Y_n,Z_n$.

The starting point in the asymptotic analysis for
$\mathbb{E}(W_n)$ is the following generating function identity:
$$
1+\sum_{n=1}^\infty\mathbb{P}(X_n\le m, Y_n\le l)q(n)x^n
=\prod_{k=1}^m\sum_{j=1}^l(1-x^{j+k-1})^{-1}, \quad l,m=1,2,... .
$$
It follows from a stronger result due to MacMahon
\cite[Section~495]{M16}. For more details and other proofs of this
result we also refer the reader to \cite[Chapter~V]{S71}. If we
keep either of the parameters $l$ and $m$ fixed, setting the other
one $:=\infty$, we obtain
\begin{eqnarray}
& & 1+\sum_{n=1}^\infty\mathbb{P}(W_n\le m)q(n)x^n
=\prod_{k=1}^m(1-x^k)^{-m}\prod_{j=m+1}^\infty
(1-x^j)^{-j} \nonumber \\
 & & =Q(x)\prod_{j=m+1}^\infty(1-x^j)^{j-m}, \quad
W_n=X_n,Y_n,Z_n. \nonumber
\end{eqnarray}
This implies the identity
\begin{equation}\label{expectlp}
\sum_{n=1}^\infty\mathbb{E}(W_n)q(n)x^n =Q(x)F_2(x),
\end{equation}
where
$$
F_2(x) =\sum_{m=0}^\infty (1-\prod_{j=m+1}^\infty (1-x^j)^{j-m}),
$$
since $\mathbb{E}(W_n)=\sum_{m=0}^{n-1} \mathbb{P}(W_n>m)$. For
the sake of convenience, we represent $F_2(x)$ in the form:
\begin{equation}\label{ftwo}
F_2(x) =\sum_{m=0}^\infty (1-e^{H_m(x)}),
\end{equation}
where
\begin{equation}\label{hm}
H_m(x) =\sum_{j>m}(j-m)\log{(1-x^j)}.
\end{equation}
Our first goal will be to find the asymptotic of $F_2(e^{-d_n})$.
Then, we will briefly sketch the verification of conditions (A)
and (B). So, in (\ref{ftwo}) we set $x=e^{-d_n}$ and break up the
sum representing $F_2(e^{-d_n})$ into three parts:
\begin{equation}\label{thr}
F_2(e^{-d_n}) =\Sigma_1+\Sigma_2+\Sigma_3,
\end{equation}
where
\begin{equation}\label{sigmaone}
\Sigma_1 =\sum_{0\le m\le N_1} (1-e^{H_m(e^{-d_n})}),
\end{equation}
\begin{equation}\label{sigmatwo}
\Sigma_2=\sum_{N_1<m\le N_2} (1-e^{H_m(e^{-d_n})}),
\end{equation}
\begin{equation}\label{sigmathree}
\Sigma_3 =\sum_{m>N_2} (1-e^{H_m(e^{-d_n})}).
\end{equation}
The choice of the the numbers $N_1$ and $N_2$ will be specified
later.

We will need first asymptotic expansions for $d_n^{-1}$,
$d_n^{-2}$ and $\log{d_n^{-2}}$. Using (\ref{dn}), it is not
difficult to show that
\begin{equation}\label{dnone}
d_n^{-1} =\left(\frac{n}{2\zeta(3)}\right)^{1/3}
+\frac{1}{36(2\zeta(3))^{2/3} n^{1/3}} +O(n^{-1/3-\beta}),
\end{equation}
\begin{equation}\label{dntwo}
d_n^{-2} =\left(\frac{n}{2\zeta(3)}\right)^{2/3}
+\frac{1}{36\zeta(3)} +O(n^{-\beta})
\end{equation}
and
\begin{equation}\label{dnlog}
\log{d_n^{-2}} =\frac{2}{3}\log{n} -\frac{2}{3}\log{(2\zeta(3))}
+O(n^{-2/3}).
\end{equation}
Next, we need to find an alternative representation for
$H_m(e^{-d_n})$. As previously, we interpret the underlying sum as
a Riemann sum with step size $d_n$. We will obtain an integral
that can be simplified using integration by parts. Thus, setting
$x=e^{-d_n}$ and $v_{m.n}=m d_n$ in (\ref{hm}), we obtain
\begin{eqnarray}\label{hmdn}
& & H_m(e^{-d_n})= d_n^{-2}\sum_{jd_n>v_{m,n}} (jd_n -v_{m,n})d_n
\log{(1-e^{-jd_n})} \nonumber \\
& & =d_n^{-2}\int_{v_{m,n}}^\infty (u-v_{m,n}) \log{(1-e^{-u})} du
+O(1) \\
& & = -\frac{d_n^{-2}}{2}\int_{v_{m,n}}^\infty
\frac{(u-v_{m,n})^2}{e^u-1} du +O(1)
=-\frac{d_n^{-2}}{2}\psi(v_{m,n}) +O(1), \quad n\to\infty,
\nonumber
\end{eqnarray}
where
$$
\psi(v)=\int_0^\infty \frac{u^2}{e^{u+v}-1} =e^{-v}\int_0^\infty
\frac{u^2 du}{e^u-e^{-v}}, \quad v\ge 0.
$$
It is easy to check, using MacLaurin formula, that
\begin{equation}\label{psi}
\psi(v) =2e^{-v} +O(ve^{-2v}), \quad v\to\infty.
\end{equation}
Furthermore, since $\log{(1-e^{-jd_n})}<0$ for all $j\ge 1$, the
sequence $\{H_m(e^{-d_n})\}_{m\ge 1}$ is monotonically increasing.
Hence. for all $m\ge N_1$,
$$
1-e^{H_m(e^{-d_n})}\ge 1-e^{H_{N_1}(e^{-d_n})}.
$$
So, using (\ref{hmdn}) and (\ref{psi}), we conclude that if
$N_1=N_1(n)\to\infty$, then the sum $\Sigma_1$ in
$(\ref{sigmaone})$ satisfies the inequalities
\begin{eqnarray}\label{none}
& & N_1\ge\Sigma_1\ge N_1-N_1 e^{H_{N_1}(e^{-d_n})} =N_1-O(N_1
 e^{-d_n^{-2} e^{-v_{N_1,n}}}) \nonumber \\
& & =N_1 -O(N_1 e^{-d_n^{-2}e^{-N_1 d_n}}).
\end{eqnarray}
The last $O$-term tends to $0$ as $n\to\infty$ if we set
\begin{equation}\label{defnone}
N_1=d_n^{-1}(\log{d_n^{-1}} -\log{\log{(N_1\log{N_1})}}).
\end{equation}
This, of course, implies that $N_1\sim
d_n^{-1}\log{d_n^{-2}}\to\infty$ (see (\ref{dnone}) and
(\ref{dnlog})). A more precise lower bound in (\ref{none}) can be
found using (\ref{defnone}) and (\ref{dn}). We have
$$
N_1 e^{-d_n^{-2}e^{-N_1 d_n}} =N_1 e^{-\log{(N_1\log{N_1})}}
=\frac{1}{\log{N_1}} =O\left(\frac{1}{\log{d_n^{-1}}}\right)
=O\left(\frac{1}{\log{n}}\right).
$$
Hence (\ref{none}) implies that $N_1\ge\Sigma_1\ge
N_1-O(1/\log{n})$, or equivalently,
\begin{equation}\label{sigmaoneas}
\Sigma_1 =N_1+O(1/\log{n}).
\end{equation}
Once the asymptotic order of $\Sigma_1$ was determined by
(\ref{sigmaoneas}), we need to find an asymptotic expression for
$N_1$ as a function of $n$. First, we analyze the
$\log{\log}$-term in (\ref{defnone}). We have
\begin{eqnarray}\label{loglog}
& & \log{\log{(N_1\log{N_1})}} =\log{(\log{N_1}+\log{\log{N_1}})}
\nonumber \\
& & =\log{\log{\left(N_1\left(1+\frac{\log{\log{N_1}}}{\log{N_1}}\right)\right)}} \nonumber \\
& &
=\log{\left(\log{N_1}+\log{\log{\left(1+\frac{\log{\log{N_1}}}{\log{N_1}}\right)}}\right)}
\nonumber \\
 & &
=\log{\left(\log{N_1}+O\left(\frac{\log{\log{N_1}}}{\log{N_1}}\right)\right)}
\nonumber \\
& & =\log{\log{N_1}}
+\log{\left(1+\frac{\log{\log{N_1}}}{\log{N_1}}\right)} \nonumber
\\
& & =\log{\log{N_1}}
+O\left(\frac{\log{\log{N_1}}}{\log{N_1}}\right).
\end{eqnarray}
Next, we will apply (\ref{dnone}) and (\ref{dntwo}). First,
(\ref{dnone}) implies the following estimate for $\log{N_1}$:
\begin{eqnarray}
& & \log{N_1} =\log{d_n^{-1}} +\log{\left(\log{d_n^{-2}}
-\log{\log{N_1}}
+O\left(\frac{\log{\log{N_1}}}{\log{N_1}}\right)\right)} \nonumber
\\
 & & =\log{\left(\left(\frac{n}{2\zeta(3)}\right)^{1/3}\right)}
 +O(n^{-1/3}) +O(\log{\log{d_n^{-2}}}) \nonumber \\
 & & =\frac{1}{3}\log{n} +O(\log{\log{n}}). \nonumber
 \end{eqnarray}
 Hence
 \begin{equation}\label{loglogtwo}
 \log{\log{N_1}} =-\log{3} +\log{\log{n}}
 +O\left(\frac{\log{\log{n}}}{\log{n}}\right).
 \end{equation}
 Combining (\ref{dnone}), (\ref{dnlog}) and (\ref{defnone}) -
 (\ref{loglogtwo}), we finally establish that
 \begin{eqnarray}\label{sigmaonefin}
 & & \Sigma_1 =\left(\left(\frac{n}{2\zeta(3)}\right)^{1/3}
 +O(n^{-1/3})\right) \nonumber \\
 & & \times\left(\frac{2}{3}\log{\left(\frac{n}{2\zeta(3)}\right)}
 +O(n^{-2/3}) +\log{3} -\log{\log{n}}
 +O\left(\frac{\log{\log{n}}}{\log{n}}\right)\right)  \\
 & & =\left(\frac{n}{2\zeta(3)}\right)^{1/3}
 \left(\frac{2}{3}\log{n} -\log{\log{n}}
  -\frac{2}{3}\log{(2\zeta(3))} +\log{3}
 +O\left(\frac{\log{\log{n}}}{\log{n}}\right)\right). \nonumber
 \end{eqnarray}
 We will estimate $\Sigma_2$ and $\Sigma_3$ (see (\ref{sigmaone})
 and (\ref{sigmatwo}), respectively) setting
 \begin{equation}\label{ntwo}
 N_2 =d_n^{-1} (\log{d_n^{-2}} +\log{\log{N_1}}).
 \end{equation}
 Using the inequality $1-e^{-u}\le u, u\ge 0$, we obtain in a
 similar way that
 \begin{eqnarray}\label{sigmatwofin}
 & & \Sigma_2\le (N_2-N_1) O(d_n^{-1}(\log{\log{N_1}})d_n^{-2}
 e^{-\log{d_n^{-2}} -\log{\log{N_1}}}) \nonumber \\
 & & =O\left(d_n^{-1}\frac{\log{\log{N_1}}}{\log{N_1}}\right)
 =O\left(n^{1/3}\frac{\log{\log{n}}}{\log{n}}\right).
 \end{eqnarray}
 Finally, by (\ref{sigmathree}) - (\ref{psi}) and (\ref{ntwo})
 \begin{eqnarray}\label{sigmathreefin}
 & & \Sigma_3\le\sum_{m>N_2} d_n^{-2} e^{-m d_n} =O(d_n^{-2}
 e^{-N_2}) \nonumber \\
 & & =O(d_n^{-2} e^{-d_n^{-1}\log{d_n^{-1}}}) =O(n^{2/3}
 e^{-2n^{1/3}(\log{n})/3}).
 \end{eqnarray}
 Now, (\ref{thr}) - (\ref{sigmathree}), (\ref{sigmaonefin}),
 (\ref{sigmatwofin}) and (\ref{sigmathreefin}) imply that
 $\Sigma_1$ presents the main contribution to the asymptotic of
 $F_2(e^{-d_n})$ and we have
 \begin{eqnarray}\label{ftwofin}
 & & F_2(e^{-d_n}) =\left(\frac{n}{2\zeta(3)}\right)^{1/3}
 \\
 & & \times\left(\frac{2}{3}\log{n} -\log{\log{n}} -\frac{2}{3}\log{(2\zeta(3))}
 +\log{3}
 +O\left(\frac{\log{\log{n}}}{\log{n}}\right)\right), \quad
 n\to\infty. \nonumber
 \end{eqnarray}
This result implies that $F_2(x)$ also satisfies condition (B). In
fact, (\ref{ftwofin}) shows that $F_2(\mid x\mid) (\ge\mid
F(x)\mid)$ is of order $O((1-\mid x\mid)^{1/3} \mid\log{(1-\mid
x\mid)\mid)} =O(e^{C/(1-\mid x\mid)^\eta})$ for any $C,\eta>0$.
The verification of condition (A) is slightly longer. It is based
on a convergence argument for $F_2(r_n)$ similar to that in
(\ref{conda}) for $F_1(r_n)$. We omit the details and refer the
reader to \cite[formulas~(3.9), (3.11)]{M06}, which imply that the
orders of growth of $H_m(x)$ and $H_m^\prime(x)$ are not larger
than the third and second powers of $(1-\mid x\mid)^{-1}$,
respectively. An argument similar to that given in the proof of
(\ref{ftwofin}) yields the required convergence. Thus, using
Theorem 1, (\ref{expectlp}) and (\ref{ftwofin}), we obtain the
following result.

\begin{proposition} If $n\to\infty$, then
$$
\mathbb{E}(W_n)\sim\kappa_2 n^{1/3}\log{n},
$$
where $W_n=X_n,Y_n,Z_n$, and
$\kappa_2=\frac{2}{3}(2\zeta(3))^{1/3} =0.4976...$.
\end{proposition}

{\it Remark 4.} In \cite{M06} is shown that all three dimensions
$X_n$, $Y_n$ and $Z_n$ of the random solid diagram with volume
$n$, appropriately normalized, converge weakly to the doubly
exponential (extreme value) distribution as $n\to\infty$.

{\it Remark 5.} It is possible to obtain more precise asymptotic
estimates (expansions) using the circle method. A kind of this
method was applied by Wright \cite{W31} who obtained an asymptotic
expansion for the numbers $q(n)$ as $n\to\infty$. His asymptotic
expansion together with a suitable expansion for $F(e^{-d_n})$
would certainly lead to better asymptotic estimates for the
expectations of various plane partition statistics.

\section*{Appendix}
\renewcommand{\theequation}{A.\arabic{equation}}
\setcounter{equation}{0}

In the Appendix we deduce Wright's formula (\ref{wright}), using
Hayman's result (\ref{qhay}).

First, by (\ref{dn}) and (\ref{bdn}) one has
\begin{equation}\label{endn}
e^{nd_n} =\exp{((2\zeta(3))^{1/3} n^{2/3} -1/36+ O(n^{-\beta}))},
\end{equation}

\begin{equation}\label{be}
\sqrt{2\pi b(e^{-d_n})} \sim\frac{(6\pi)^{1/2} n^{2/3}}
{(2\zeta(3))^{1/6}}.
\end{equation}

An asymptotic expression for $Q(e^{-d_n})$ can be obtained using a
general lemma due to Meinardus \cite{M54} (see also
\cite[Lemma~6.1]{A76}). Since the Dirichlet generating series for
the plane partitions is $\zeta(z-1)$, we get
\begin{eqnarray}
& & Q(e^{-d_n}) =\exp{(\zeta(3)d_n^{-2} -\zeta(-1)\log{d_n}
+\zeta^\prime(-1) +O(d_n^{\beta_1}))} \nonumber \\
& & =\exp{(\zeta(3)d_n^{-2} +\frac{1}{12}\log{d_n} +2\gamma
+O(d_n^{\beta_1}))} \nonumber
\end{eqnarray}
where $0<\beta_1<1$ and $\gamma$ is given by (\ref{glakin}) (more
details on the values of $\zeta(-1)$ and $\zeta^\prime(-1)$ can be
found in \cite[Section~13.13]{WW27} and \cite[Section~2.15]{F03}).
Using (\ref{dntwo}) and (\ref{dnlog}), after some algebraic
manipulations, we obtain
\begin{equation}\label{qedn}
Q(e^{-d_n}) =\left(\frac{2\zeta(3)}{n}\right)^{1/36}
\exp{((\zeta(3))^{1/3}(n/2)^{2/3} +1/36 +O(n^{-\beta_1}))}.
\end{equation}
Combining (\ref{endn}) - (\ref{qedn}), we find that
\begin{eqnarray}
& & q(n) \sim\left(\frac{2\zeta(3)}{n}\right)^{1/36}
\frac{\exp{((2\zeta(3))^{1/3}-1/36+(\zeta(3))^{1/3}(n/2)^{2/3}+1/36+2\gamma)}}{(3\pi)^{1/2} n^{2/3}/(2\zeta(3))^{1/6}} \nonumber \\
& & =\frac{(\zeta(3))^{1/6+1/36} n^{-1/36-2/3}}
{2^{1/2-1/6-1/36}(3\pi)^{1/2}}
\exp{(3(\zeta(3))^{1/3}(n/2)^{2/3}+2\gamma)} \nonumber \\
& & =\frac{(\zeta(3))^{7/36}}{2^{11/36}(3\pi)^{1/2}} n^{-25/36}
\exp{(3(\zeta(3))^{1/3}(n/2)^{2/3}+2\gamma)}.\quad \Box \nonumber
\end{eqnarray}

\end{document}